\documentclass[]{article}

\usepackage{amsmath}
\usepackage{amssymb}
\usepackage[ansinew]{inputenc}

\def\Cset{\mathbb{C}}
\def\Nset{\mathbb{N}}

\def\Rset{\mathbb{R}}

\newtheorem{theorem}{Theorem}
\newtheorem{lemma}{Lemma}
\newtheorem{example}{Example}
\newtheorem{cor}{Corollary}
\newtheorem{defi}{Definition}
\newtheorem{remark}{Remark}
\newtheorem{prop}{Proposition}
\newcommand{\proof}{\bfseries Proof: $\quad$\mdseries}
\newcommand{\finishproof}{\begin{flushright} $\Box$ \end{flushright}}

\date{\today}

\begin{document}

\title{$L^{2}$-spectral gaps, weak-reversible and very weak-reversible Markov chains}
\author{Achim Wübker\thanks{e-mail: awuebker@mathematik.uni-osnabrueck.de}\\Institute of Mathematics\\Albrechtstraße 28 a, 49076 Osnabrück \and Zakhar Kabluchko\thanks{e-mail: kabluch@math.uni-goettingen.de}\\Institute of Mathematical Stochastic\\Goldschmidtstr. 7, 37077 Göttingen}

\maketitle

\begin{abstract}
The theory of $L^2$-spectral gaps for reversible Markov chains has been studied by many authors. In
this paper we consider positive recurrent general state space Markov chains with stationary transition
probabilities. Replacing the assumption of reversibility by a less strong one, we still obtain a simple necessary and sufficient condition for the spectral gap property of the associated Markov operator in
terms of isoperimetric constant. Moreover, we define a new sequence of isoperimetric constants which provides a necessary and sufficient condition for the existence of a spectral gap in a very general setting. Finally, these results are used to obtain simple sufficient
conditions for the existence of a spectral gap in terms of the first and second order transition probabilities.
\end{abstract}

\section{Introduction}
Let $\xi_1,\xi_2,\ldots$ be a time discrete and time homogeneous positive recurrent Markov chain on an arbitrary state space $(\Omega,\mathcal{F})$ with transition kernel $p(\cdot,\cdot)$ and 
uniquely determined invariant measure $\pi$. The main question addressed is the existence of an $L^2(\pi)$-spectral gap. We are interested in conditions that ensure  

\[\lim_{n\rightarrow\infty}\sup_{f\in L_{0,1}^{2}(\pi)}||P^{n}f||_{2}^{\frac{1}{n}}<1,\]
where 
\[
Pf(x):=\int_{\Omega}f(y)P(x,dy),\,\,f\in L^{2}(\pi)
\]
and $L_{0,1}^{2}(\pi):=\{f\in L^{2}(\pi):\int_{\Omega}f(x)\pi(dx)=0,\,\int_{\Omega}f(x)^{2}\pi(dx)=1\}$.\\
For reversible Markov chains, simple conditions equivalent to the spectral gap property
are known in the time discrete case \cite{wuebker} as well as in the time continuous case \cite{chen}, \cite{chenbuch}, \cite{chenwang}, \cite{diaconis3}, \cite{diaconis4}, \cite{lawler}. These conditions are given in terms of 
isoperimetric constants \cite{chen}, \cite{chenbuch}, \cite{chenwang}, \cite{dodziuk}, \cite {lawler}, \cite{wuebker} or in terms of geometric constants \cite{diaconis3}.
So far, without the assumption of reversibility, only a few results under very strong 
conditions are known \cite{fill}. There are mainly two reasons for this: First, non 
reversibility implies that Dirichlet form techniques can only be
successfully applied by using certain reversibilization procedure (e.g. \cite{lawler}, \cite{diaconis1}, \cite{fill}).  
But then one only obtains information about the real part of the spectrum.
The second problem is that there is no canonical generalization procedure. Since reversibility corresponds
to the self-adjointness of the associated Markov operator $P$, one first might try to generalize
the results obtained in \cite{wuebker} to normal Markov operators. We will show by example that in 
general this approach won't be successful. In order to tackle the problem, we compare the 
isoperimetric constant associated to $P^{2}$ with these associated to $P^{\ast}P$ and $PP^{\ast}$ (where $P^{\ast}$ denotes the adjoint of $P$). This technique turns out to be appropriate to study spectral gap properties of $P$ in the non-reversible case if the Markov chain satisfies a condition that is much weaker than reversibility, but is at least related to the latter and is called very weak reversibility for that reason. In fact, this method enables us to prove Theorem \ref{three}, which generalizes the results obtained in \cite{wuebker} and can be seen as the core of this paper. Moreover it is possible to use the theory developed in this 
paper to obtain estimates for spectral gaps of non-reversible Markov chains.


Let us introduce the basic notations and recall some known facts. Unless stated otherwise, we consider a positive recurrent time homogeneous Markov chain $\xi_1,\xi_2,\ldots$ with arbitrary state space $(\Omega,\mathcal{F})$, transition kernel $p(x,dy)$ and uniquely determined invariant probability measure $\pi$.
We say that the chain $\xi_1,\xi_2,\ldots$ is reversible, if for all $A,B\in\mathcal{F}$
\[
Q(A,B):=\int_{A}\pi(dx)\int_{B}p(x,dy)=\int_{B}\pi(dx)\int_{A}p(x,dy)=Q(B,A)=:\tilde{Q}(A,B),
\]
so $Q$, $\tilde{Q}$ are measures on $\Omega^{2}$.
Alternatively, reversibility can be stated as
\[\frac{dQ}{d\tilde{Q}}(x,y)=1,\,\,\,\tilde{Q}\,\,a.s.\]
In order to give a natural generalization of this definition, let us, for simplicity, assume that
$Q$ and $\tilde{Q}$ are equivalent measures in the Radon-Nikodym-sense.   
\begin{defi}\label{weakrev}
We say the Markov chain $\xi_1,\xi_2,\ldots$ is weak reversible of order $n\in\Nset$ if there exists $C\in[1,\infty)$ such that 
\[\frac{1}{C}\le\frac{dQ^{(n)}}{d\tilde{Q}^{(n)}}\le C,\,\,\,\tilde{Q}^{(n)}\,\,a.s.,\]
where $Q^{(n)}(A,B):=\int_{A}\pi(dx)p^{n}(x,B)$ and $\tilde{Q}^{(n)}(A,B)=Q^{(n)}(B,A)$.
\end{defi}

For $C=n=1$ this is exactly the definition of reversibility. \\

The definition above can be generalized to
\begin{defi}
We say the Markov chain $\xi_1,\xi_2,\ldots$ is very weak reversible of order $n\in\Nset$ if there exists $q\in(1,\infty]$ such that
\begin{equation}
\textrm{esssup}_{y\in\Omega}||\frac{dQ^{(n)}}{d\tilde{Q^{(n)}}}(\cdot,y)||_{L^{q}(p(y,\cdot))}<\infty\,\,\,\pi\,a.s.
\end{equation}
\end{defi}
From the definition it is immediately clear that weak reversible Markov chains are very weak reversible,
but in general the converse is not true.

In the sequel we will need the following families of isoperimetric constants: 
\[\label{isoneu}
k_{n}:=\inf_{A\in\mathcal{F}}k_{n}(A),\quad k_{n}(A):=\frac{1}{\pi(A)\pi(A^{c})}\int_{A}p^{n}(x,A^{c})\pi(dx), \,\,n\in\Nset,
\]
and
\[k_{P^{\ast^{n}}P^{n}}:=\inf_{A\in\mathcal{F}}k_{P^{\ast^{n}}P^{n}}(A):=\inf_{A\in\mathcal{F}}\frac{1}{\pi(A)\pi(A^{c})}\int_{A}P^{\ast^ {n}}P^{n}1_{A^{c}}(x)\pi(dx),\]
where $P^{\ast}$ is the adjoint operator of $P$ considered on $L^{2}(\pi)$.\\
In \cite{wuebker} we saw that aperiodicity of the Markov chain $\xi_1,\xi_{n+1},\ldots$ can be measured by the constants
\[K_n:=\sup_{A\in\mathcal{F}}k_{n}(A), n\in\Nset.\]
We call an operator $P$ positive if $Pf\ge 0$ for all $f\ge 0$.
One can show that there exists a measure $\mu$ on $\Omega^{2}$ such that for all $f,g\in L^{2}(\pi)$
\begin{equation}\label{desi}
<f,P^{\ast}Pg>_{\pi}:=\int_{\Omega}f(x)P^{\ast}Pg(x)\pi(dx)
=\int_{\Omega}\int_{\Omega}\mu(dx,dy)g(x)f(y).
\end{equation}
From this, we obtain that
\[
<f,(Id-P^{\ast}P)g>_{\pi}=\frac{1}{2}\int_{\Omega}\int_{\Omega}(g(y)-g(x))(f(y)-f(x))\mu(dx,dy)
\]
and therefore especially
\begin{equation}\label{deshalbgehtss}
<f,(Id-P^{\ast}P)f>_{\pi}=\frac{1}{2}\int_{\Omega}\int_{\Omega}(f(y)-f(x))^{2}\mu(dx,dy).
\end{equation}
This representation will be used in the proof of Theorem \ref{fife}.\\
If we replace $P^{\ast}P$ by an arbitrary positive and self adjoint operator $A$, a 
representation as in (\ref{deshalbgehtss}) holds true, provided that the according probability
space $(\Omega,\mathcal{F},\pi)$ satisfies some weak conditions (see \cite{fukushima}).\\
The spectrum $\sigma(P)$ of $P$ may be divided into three disjoint components,
\begin{equation}
\sigma(P)=\sigma_{p}(P)\cup\sigma_{c}(P)\cup\sigma_{r}(P),
\end{equation}
where $\sigma_{p}(P)$ is called the eigenspectrum (or discrete spectrum), $\sigma_{c}(P)$ the continuous
spectrum and $\sigma_{r}(P)$ the residual spectrum. The approximate spectrum $\sigma_{ap}(P)$, which
is defined to be the collection of complex numbers which satisfy
\begin{equation}
\lim_{n\rightarrow\infty}||(P-\lambda Id)f_{n}||\rightarrow 0
\end{equation}
for some sequence $(f_{n})$ with $||f_{n}||=1$ for all $n\in\Nset$.
It is well know that 
\begin{equation}
\sigma_{p}(P)\cup\sigma_{c}(P)\subset\sigma_{ap}(P)\subset\sigma(P)
\end{equation}
and in the case where $P$ is normal (i.e. $P^{\ast}P=PP^{\ast}$), it holds (see \cite{birman})
\begin{equation}
\sigma_{p}(P)\cup\sigma_{c}(P)=\sigma_{ap}(P)=\sigma(P).
\end{equation}


\begin{remark}
The definitions of $k_n$ and $K_n$ are taken from \cite{wuebker}, that of weak reversible, very weak
reversible and
$k_{P^{\ast^{n}}P^{n}}$ appears to be new. 
\end{remark}

\section{Spectral theory for general and weak reversible Markov chains}
Let us start with the following, probably well-known Proposition:

\begin{prop}\label{one}
Let $\xi_1,\xi_2,\ldots$ be a positive recurrent Markov chain. Then the following two statements are equivalent: 
\begin{enumerate}
\item
$P$ has an $L^{2}(\pi)$-spectral gap.
\item
\begin{equation}\label{hau5}
\exists n_{0}\in\Nset:\,\,k_{P^{\ast^{n_{0}}}P^{n_{0}}}>0.
\end{equation}
\end{enumerate}
If (\ref{hau5}) is satisfied, we obtain the following estimate for the spectral radius $r$ of $P$ on $L_{0,1}^{2}(\pi)$:
\begin{equation}\label{charspeckadj}
\sigma(P)\subset B_{r}(0):=\{x\in\Cset:||x||_{2}\le r\},\, r=\left(\sqrt{1-\frac{\kappa}{8}k_{P^{\ast^{n_{0}}}P^{n_{0}}}^{2}}\right)^{\frac{1}{n_0}},
\end{equation}
$n_{0}$ as in (\ref{hau5})
\end{prop}
At the end of the paper we provide an example where this can be used to decide whether the operator
$P$ has an $L^{2}(\pi)$-spectral gap or not. \\
The next theorem provides a necessary and sufficient condition for the existence of an $L^{2}(\pi)$-spectral gap. The interesting fact is that only one isoperimetric constant is needed in order to
establish the existence of a spectral gap: 

\begin{theorem}\label{two}
Let us assume that $\xi_1,\xi_2,\ldots$ is very weak reversible of order $n$. Then the following two
conditions are equivalent:
\begin{enumerate}
\item
$P$ has a $L^2(\pi)$-spectral gap.
\item
$k_{2n}>0.$
\end{enumerate}
\end{theorem}
Now let us state the main theorem of this paper:
\begin{theorem}\label{three}
Let us assume that $\xi_1,\xi_2,\ldots$ is weak reversible of order $n$. Then the following three conditions are equivalent:

\begin{enumerate}
\item
$P$ has a $L^2(\pi)$-spectral gap.
\item
$k_{2n}>0.$
\item
$0<k_n\le K_n<2.$
\end{enumerate}
\end{theorem}
We see here that under the assumption of weak reversibility the isoperimetric constants $k_n$ and $K_n$ are still an appropriate tool for analyzing the existence and the size of a spectral gap. This shows that the operator $P$ may be far away from being reversible, but the concept of isoperimetric constants
can still be applied in order to obtain estimate for the spectral gap.
Note that for $n=1$, Theorem \ref{three} improves the result in \cite{wuebker}, where reversibility is assumed. \\
In order to state the next theorem, we need to introduce the following notations:
Let $A_n\in\mathcal{F}$ be such that $\lim_{n\rightarrow\infty}k_2(A_n)=k_2$. Without loss of generality we may assume that
$\pi(A_n)\le\frac{1}{2}$. Let $\pi_{A}(B):=\frac{\pi(A\cap B)}{\pi(A)}$ and $f_{A_n}(x):=\frac{p(x,A_n^{c})}{p^{2}(x,A_n^{c})}$.
Assume that there exists $p\in (1,\infty]$ and a sequence $A_n$ as above such that 
\begin{equation}\label{theo4}
\sup_{n\in\Nset}||1_{A_n}f_{A_n}||_{L^{p}(\pi_{A_n})}<\infty.
\end{equation}
\begin{theorem}\label{four}
Let $\xi_1,\xi_2,\ldots$ be a very weak reversible Markov chain such that (\ref{theo4}) holds. Then $\xi_1,\xi_2,\ldots$ has a spectral gap if and only if
\begin{equation}
k:=k_1>0.
\end{equation}
\end{theorem}
Now we should have a closer look at the assumption (\ref{theo4}). Let us consider an $\epsilon$-lazy
Markov chain, i.e. a Markov chain that satisfies
\begin{equation}\label{mcmc}
p(x,x)\ge\epsilon
\end{equation}
for an $\epsilon>0$ and all $x\in \Omega$. Then we have that $f_{A_n}\le\frac{1}{\epsilon}$ and so
(\ref{theo4})is trivially fulfilled. But in fact, for MCMC simulations condition (\ref{mcmc}) is often satisfied or the chain can be defined in such a way that
(\ref{mcmc}) is satisfied (see e.g. \cite{bremaud}). For example, it is often assumed that
$p(x,x)\ge\frac{1}{2}$ in order to have $\mathcal{R}(\sigma(P))\ge 0$, i.e. all values in the spectrum of $P$ should have a non negative real part.
Another advantage of the Theorem \ref{four} is
that it is possible to give lower bounds for the size of the spectral gap by combining it with Proposition \ref{one}, Lemma \ref{iteriertintegral} and Lemma \ref{nele}.   \\

Let us state the next theorem, which is closely connected to Theorems 2.1. and 2.3. in \cite{lawler}.

\begin{theorem}\label{fife}
Let $\xi_1,\xi_2,\ldots$ be a positive recurrent Markov chain such that for $P$ we have that
$\sigma_{ap}(P)=\sigma(P)$.  Moreover, let us assume that
\begin{equation}\label{posii}
P+P^{\ast}-P^{\ast}P
\end{equation}
is a positive operator. Then, for the spectrum $\sigma(P)$ of $P$ acting on $L_{0,1}^{2}(\pi)$, we obtain that 
\begin{equation}
\sigma(P)\bigcap B_{\sqrt{\frac{\kappa}{8}}k^2}(1)=\emptyset.
\end{equation}
where $\kappa$ denotes the constant
introduced in \cite{lawler}.
If in addition $P$ is assumed to be self adjoint, we obtain
\[\sigma(P)\subset[0,1-\sqrt{\frac{\kappa}{8}}k^{2}].\] 
\end{theorem}
Here it is interesting that we obtain a sufficient condition for the existence of a spectral gap
in terms of $k$ without using any reversibility assumptions on $p$.

Under certain conditions it is possible to compare this with the results of Lawler and Sokal \cite{lawler}. To this end let us
remark that the result above with slightly modifications can be proved for positive recurrent continuous-time Markovian jump processes
in the same way replacing $p(x,dy)$ by
$j(x,dy)1_{x\not=y}$, $(Id-P)f$ by $\tilde{J}f:=\int_{\Omega}j(x,dy)1_{x\not=y}(f(x)-f(y))$, (\ref{posii}) by \\ $J+J^{\ast}-J^{\ast}J$ to be positive and $k²$ by $k_{J+J^{\ast}-J^{\ast}J}$. Here $j$ is
the transition rate function corresponding to the jump process. In the continuous time case theorem 2.3. in Lawler and Sokal \cite{lawler} states that for reversible Markov chains it holds
\[\mathcal{R}(\sigma(\tilde{J}))\ge \frac{\kappa}{8M}k^2.\]
Our result states that under the additional assumptions of $\sigma(P)=\sigma_{ap}(P)$ and positivity of $J+J^{\ast}-J^{\ast}J$ one can skip the reversibility assumption to obtain
\[\sigma(\tilde{J})\subset B_{\sqrt{\frac{\kappa}{8M}}k_{J+J^{\ast}-J^{\ast}J}}(0)^{c}.\]
This implies that the size of the gap decreases asymptotically with rate $\sqrt{\frac{1}{M}}$ for $M\rightarrow\infty$, which 
is substantially slower compared to $\frac{1}{M}$ and therefore improves the results obtained in \cite{chenwang} and \cite{lawler} for reversible Markov chains in that direction, provided $k_{J+J^{\ast}-J^{\ast}J}>0$ .\\
For reversible Markov chains, equation (\ref{posii}) can be expressed in terms of the transition probabilities and is equivalent to
\begin{equation}\label{verystronga}
2p(x,\cdot)-p^{2}(x,\cdot)\ge 0 \,\,\pi \,\,a.s.
\end{equation}
Moreover, from the spectral mapping theorem it follows that (\ref{posii}) implies
\[\sigma(P)\subset[1-\sqrt{2},1].\]
These considerations provide the framework where Theorem \ref{fife} can be applied. Using (\ref{verystronga}) it is easy to see that for reversible chains with countable state space the condition (\ref{posii}) can only be satisfied for Markov chain which enters from every starting point every other point of the state space with positive probability.


The following observation will show why it is interesting to consider weak (very weak) reversible Markov chains of higher orders:
\begin{equation}\label{gegenbeispiel2}
P=\left(\begin{array}{ccc}
0&1&0\\
0&0&1\\
\frac{1}{2}&0&\frac{1}{2}\\
\end{array}
\right).
\end{equation}
The invariant probability measure $\pi$ is given by $\pi=(\frac{1}{4},\frac{1}{4},\frac{1}{2})$.
One can check that 
\[
P^{\ast}P=\left(\begin{array}{ccc}
\frac{1}{2}&0&\frac{1}{2}\\
0&1&0\\
\frac{1}{4}&0&\frac{3}{4}\\
\end{array}
\right).
\]
Since $P^{\ast}P$ is not ergodic, $k_{P^{\ast}P}=0$ and hence it follows from Lemma $4$ of \cite{wuebker}
that $k_{(P^{\ast}P)^{n}}=0\,\,\forall n\in\Nset$. But of course, this chain has the spectral gap property, which will be obvious by considering $k_{P^{\ast}²P²}$  \\
The problem is that without any restrictions on the transition probabilities, deterministic set movements of order larger than 2 are possible. We will see that the weak reversibility property excludes such behavior and will facilitate a comparison between $k_{2}$ and $k_{P^{\ast}P},k_{PP^{\ast}}$.\\



\section{Proofs of the theorems}\label{}
Let us start with the proof of Proposition \ref{one}.
\proof
Since $P^{\ast^{n}}P^{n}$ is positive and self adjoint and since (\ref{desi}) is satisfied, we can use the proof of Theorem 2.1 in \cite{lawler}. Therefore we obtain for all $f\in L^{2}_{0,1}(\pi)$:
\[\frac{\kappa}{8}k_{P^{\ast^{n}}P^{n}}^{2}\le\inf_{f\in L^{2}_{0,1}(\pi)}<f,f-P^{\ast^{n}}P^{n}f>_{2}\le k_{P^{\ast^{n}}P^{n}}.\]
This is equivalent to
\begin{equation}\label{zweimal1}
\sqrt{1-\frac{\kappa}{8}k_{P^{\ast^{n}}P^{n}}^{2}}\ge ||P^{n}f||_{2}\ge \sqrt {1-k_{P^{\ast^{n}}P^{n}}},
\end{equation}
It is not difficult to see that for all $n$ in $\Nset$ we have that $1-k_{P^{\ast^{n}}P^{n}}\ge 0$, so the right hand side of (\ref{zweimal1}) is well defined. But from this inequality the necessity of (\ref{hau5}) follows. \\
On the other hand, if (\ref{hau5}) is fulfilled for some $n_0\in\Nset$, we obtain by using the left hand side of inequality (\ref{zweimal1})

\[||P^{n_{0}n}||_{2}\le ||P^{n_{0}}||_{2}^{n}\le \left(\sqrt{1-\frac{\kappa}{8}k_{P^{\ast^{n_{0}}}P^{n_{0}}}^{2}}\right)^n,\]
Having in mind that $||P^n||_2$ is monotonic decreasing in $n$ the desired estimate follows by taking the $n_{0}n$-th root. 
\finishproof

\begin{lemma}\label{iteriertintegral}
We consider the Markov chain $\xi_1,\xi_2,\ldots$. The following inequalities hold:
\begin{equation}
k_{P^{\ast^{n_{0}}}P^{n_{0}}}\le 2^{1/p}\textrm{esssup}_{y\in\Omega}||\frac{dQ}{d\tilde{Q}}(\cdot,y)||_{L^{p}(p(y,\cdot))}k_{2n_{0}}^{1/q}\end{equation}
\begin{equation}
k_{2n_{0}}\le 2^{1/p}\textrm{esssup}_{y\in\Omega}||\frac{dQ}{d\tilde{Q}}(\cdot,y)||_{L^p(p(y,\cdot))}k_{P^{n_{0}}P^{\ast^{n_{0}}}}^{1/q}.
\end{equation}
\end{lemma}
\proof
Without loss of generality we can assume that $\pi(A^{c})\ge\frac{1}{2}$.
The first inequality can be seen as follows:
\begin{eqnarray}
\pi(A)\pi(A^{c})k_{P^{\ast}P}(A)&=&\int_{A}P^{\ast}P1_{A^{c}}(x)\pi(dx)=\int_{\Omega}p(x,A)p(x,A^{c})\pi(dx)\nonumber\\
&=&\pi(A)\int_{\Omega}\int_{A}\frac{dQ}{d\tilde{Q}}(x,y)p(x,A^{c})\frac{\pi(dy)}{\pi(A)}p(y,dx)\nonumber\\
&\le& \pi(A)\left(\int_{\Omega}\int_{A}\frac{dQ}{d\tilde{Q}}(x,y)^{p}\frac{\pi(dy)}{\pi(A)}p(y,dx) \right)^{\frac{1}{p}}\nonumber\\
&&\qquad\qquad\qquad\cdot\left(\int_{\Omega}\int_{A}p(x,A^{c})^{q}\frac{\pi(dy)}{\pi(A)}p(y,dx) \right)^{\frac{1}{q}}\nonumber\\
&\le&\pi(A)\left(\int_{A}||\frac{dQ}{d\tilde{Q}}(\cdot,y)||_{L^p(p(y,\cdot))}^{p}\frac{\pi(dy)}{\pi(A)}\right)^{\frac{1}{p}} \left(\pi(A^{c})k_{2}(A)\right)^{\frac{1}{q}}\nonumber\\
&\le&\pi(A)\pi(A^{c})2^{\frac{1}{p}}\textrm{esssup}_{y\in\Omega}||\frac{dQ}{d\tilde{Q}}(\cdot,y)||_{L^p(p(y,\cdot))}k_{2}(A)^{\frac{1}{q}}.\nonumber\\
\end{eqnarray}
Let us turn over to second inequality:
\begin{eqnarray}
\pi(A)\pi(A^{c})k_{2}(A)&=&\pi(A)\pi(A^{c})k_{2}(A^{c})=\int_{\Omega}P^{\ast}1_{A^{c}}(x)P1_{A}(x)\pi(dx)\nonumber\\
&=&\pi(A)\int_{\Omega}\int_{A}P^{\ast}1_{A^{c}}(x)p(x,dy)\frac{\pi(dx)}{\pi(A)}\nonumber\\
&\le& \pi(A)\textrm{esssup}_{y\in\Omega}||\frac{dQ}{d\tilde{Q}}(\cdot,y)||_{L^p(p(y,\cdot))}\left(\int_{A}\int_{\Omega}P^{\ast}1_{A^{c}}(x)^{q}p(y,dx)\frac{\pi(dy)}{\pi(A)}\right)^{\frac{1}{q}}\nonumber\\
&\le& \pi(A)\pi(A^{c})2^{\frac{1}{p}}\textrm{esssup}_{y\in\Omega}||\frac{dQ}{d\tilde{Q}}(\cdot,y)||_{L^p(p(y,\cdot))}k_{PP^{\ast}}(A)^{\frac{1}{q}}.\nonumber\\
\end{eqnarray}
\finishproof

From here we start the proof of Theorem \ref{two}. That the existence of an $n_{0}$ with $k_{2n_{0}}>0$ is
necessary was shown in \cite{wuebker}. Assume now that there exists $n_{0}$ such that $k_{2n_{0}}>0$. Then we have by Lemma \ref{iteriertintegral} and weak reversibility that $k_{P^{n_{0}}P^{\ast^{n_{0}}}}>0$. From Theorem \ref{one} we obtain that $P^{\ast}$ has a spectral gap.
But it is well-known that $\sigma(P)=\sigma(P^{\ast})$. So this implies that $P$ has an $L^{2}(\pi)$-spectral gap. 
\finishproof

\begin{cor}\label{corone}
Assume that $\xi_1,\xi_2,\ldots$ is a weak reversible Markov chain of order $n_0$ with reversibility constant $C$. Then we obtain the following estimates:
\begin{equation}
k_{P^{\ast^{n_{0}}P^{n_{0}}}}\le C k_{2n_{0}},\,\,\,k_{2n_0}\le  C k_{P^{n_{0}}P^{\ast^{n_{0}}}}
\end{equation}
\end{cor}
\proof
Use Lemma \ref{iteriertintegral} and choose $p=\infty$, $q=1$.
\finishproof

\begin{cor}
Let us assume that $\xi_1,\xi_2,\ldots$ is a weak reversible Markov chain and that the transition kernel of the reversed 
Markov chain $\xi_1^{\ast},\xi_2^{\ast},\ldots$ is given by $p^{\ast}(\cdot,\cdot)$.
Then we have the following inequalities:
\begin{equation}
\frac{k_{2n_{0}}^{q}}{2^{q/p}\textrm{esssup}_{y\in\Omega}||\frac{d\tilde{Q}}{dQ}(\cdot,x)||_{p,p(x,\cdot)}^{q}}\ge k_{P^{\ast^{n_{0}}}P^{n_{0}}}\ge 2^{1/p}\textrm{esssup}_{y\in\Omega}||\frac{dQ}{d\tilde{Q}}(\cdot,x)||_{p,p(x,\cdot)}k_{2n_{0}}^{1/q}
\end{equation}
\end{cor}
\proof
The first inequality is obtained by changing the roles of $P$ and $P^{\ast}$ in the proof of Lemma \ref{iteriertintegral}. 
\finishproof

For proving Theorem \ref{three}, we have to establish the following important Lemma: 

\begin{lemma}\label{main}
Let $\xi_1, \xi_2, \ldots$ be a stationary and weak reversible MC of order $n$ on an arbitrary state space $(\Omega,\mathcal{F},\pi)$. Let $C_R$ be the reversibility constant associated to the
MC. Then we obtain the following estimate for $k_{2n}$:
\begin{eqnarray}\label{jakob7}
k_{2n}&\ge&\sup_{\delta,\epsilon_{1},\epsilon_{2},\epsilon\in\mathbb{R}_{+}}\min \left[\frac{k_{n}^{2}}{16}\delta,\frac{k_{n}}{4}(
\epsilon_{1}\epsilon_{2}(1-\delta)-C_{R}\delta),\right.\nonumber\\
&&\,\,\,\,\,\,\,\,\left.\left(k_{n}\left(\frac{(2-\epsilon)(1-\epsilon_1)(1-\epsilon_2)(1-\delta)}{(1-\epsilon)K_n}-\frac{1}{1-\epsilon}\right)-\frac{\epsilon}{1-\epsilon}\right)\epsilon\right].
\end{eqnarray}
\end{lemma}

The Lemma \ref{main} is closely related to Lemma 6 of \cite{wuebker} in the way that
the assumption of reversibility is replaced by weak reversibility.
Fortunately, the proof given in \cite{wuebker} is somehow stable under the milder condition of
weak reversibility. This abbreviates the proof and we present only the part where weak 
reversibility is needed.  

\proof
Without loss of generality let us assume that $n=1$ (otherwise argue with $p^{n}(\cdot,\cdot)$ instead of
$p(\cdot,\cdot)$). As in \cite{wuebker} one can show that without loss of generality we can choose $A\in\mathcal{F}$ such that $\pi(A)\le\frac{1}{2}$. Let us consider the following sets:
\[A_{\frac{k}{4}}:=\{y\in A:p(y,A^{c})\ge\frac{k}{4}\},\,\,C:=A_{\frac{k}{4}}^{c}\cap A,\]
\[B_{\epsilon_1}:=\{x\in A_{\frac{k}{4}}:p(x,A^{c})<1-\epsilon_1\}, \]
As in \cite{wuebker} we distinguish three cases. The estimates obtained there for the first and the third case are valid without the assumption of reversibility. So let us remind the estimates obtained in \cite{wuebker}:
\begin{enumerate}
\item[Case 1:]
\[k_{2}(A)\ge\frac{k^{2}}{16}\delta_{A} \,\,\,\mbox{for }\pi(C)\ge\delta_{A}\pi(A).\]
\item[Case 3:]
\begin{eqnarray}
k_{2}(A)&\ge&\epsilon\left(k\left(\frac{(2-\epsilon)(1-\epsilon_1)(1-\epsilon_2)(1-\delta_{A})}{(1-\epsilon)K}-\frac{1}{1-\epsilon}\right)-\frac{\epsilon}{1-\epsilon}\right)\nonumber\\
&&\mbox{for }\pi(C)\le\delta_{A}\pi(A),\,\,\,\, \pi(B_{\epsilon_1})\le\epsilon_2\pi(A_{\frac{k}{4}}).\nonumber\\
\end{eqnarray}
\end{enumerate}
The second case differs from that in \cite{wuebker}, so let us assume that
$\pi(C)\le\delta_{A}\pi(A)$ and an $\epsilon_2> 0$ exists 
such that 
\begin{equation}\label{nachge}
\pi(B_{\epsilon_1})\ge\epsilon_2\pi(A_{\frac{k}{4}}).
\end{equation}
With these assumptions it follows that
\begin{eqnarray}\label{jakob3}
k_2(A)&=&\frac{1}{\pi(A)\pi(A^{c})}\int_{A}\pi(dx)p^{2}(x,A^{c})\ge\frac{1}{\pi(A)\pi(A^{c})}\int_{B_{\epsilon_1}}\pi(dx)p^{2}(x,A^{c})\nonumber\\
&\ge&\frac{1}{\pi(A)\pi(A^{c})}\int_{B_{\epsilon_1}}\pi(dx)\int_{A_{\frac{k}{4}}}p(x,dy)p(y,A^{c})\nonumber\\
&\ge&\frac{k}{4}\frac{1}{\pi(A)\pi(A^{c})}\int_{A_{\frac{k}{4}}}\pi(dx)p(x,B_{\epsilon_1})\ge\frac{k}{4C_{R}}\frac{1}{\pi(A)\pi(A^{c})}\int_{\epsilon_1}\pi(dx)p(x,B_{A_{\frac{k}{4}}}).\nonumber\\
\end{eqnarray}
Moreover, we have
\begin{enumerate}
\item
\[
\int_{B_{\epsilon_1}}\pi(dx)p(x,C)\le C_{R}\int_{C}\pi(dx)p(x,{B_{\epsilon_1}})\le C_{R} \pi(C)\le C_{R}\delta_{A}\pi(A).
\]
\item
\[
\int_{B_{\epsilon_1}}\pi(dx)p(x,A)\ge\epsilon_1\pi(B_{\epsilon_1})\ge\epsilon_1\epsilon_2\pi(A_{\frac{k}{4}})\ge\epsilon_1\epsilon_2 (1-\delta_{A})\pi(A).
\]
\end{enumerate}
Subtracting the first inequality from the second we obtain
\[\int_{B_{\epsilon_1}}\pi(dx)p(x,A_{\frac{k}{4}})\ge (\epsilon_2\epsilon_1(1-\delta_{A})-C_{R}\delta_{A})\pi(A).\]
This inserted into (\ref{jakob3}) yields
\begin{equation}
k_{2}(A)\ge\frac{k}{4 C_{R}}(\epsilon_2\epsilon_1(1-\delta_{A})-C_{R}\delta_{A}).
\end{equation}
Now gluing the three cases together yields
\begin{eqnarray}
k_2&\ge&\sup_{\delta,\epsilon_{1},\epsilon_{2},\epsilon\in\Rset_{+}}\min \left[\frac{k^{2}}{16}\delta,\frac{k}{4}(
\epsilon_{1}\epsilon_{2}(1-\delta)-C_{R}\delta),\right.\nonumber\\
&&\,\,\,\,\,\,\,\,\left.\left(k\left(\frac{(2-\epsilon)(1-\epsilon_1)(1-\epsilon_2)(1-\delta)}{(1-\epsilon)K}-\frac{1}{1-\epsilon}\right)-\frac{\epsilon}{1-\epsilon}\right)\epsilon\right].
\end{eqnarray}
This proves Lemma \ref{main}.
\finishproof

Theorem \ref{three} follows now from Theorem \ref{two} and 
Lemma \ref{main}, since the right hand side of (\ref{jakob7}) can be bounded from below if $0<k_n\le K_n<2$ (see \cite{wuebker}). \\

In order to prove Theorem \ref{four}, we should establish the following lemma.
\begin{lemma}\label{nele}
For all $p\in[1,\infty]$ and 
$q$ such that $\frac{1}{p}+\frac{1}{q}=1$ we have the following inequality:
\begin{equation}\label{nochungl}
k_{2}\ge (\frac{1}{2})^{\frac{q}{p}}\frac{k^{q}}{\sup_{n\in\Nset}||1_{A_n}f_{A_n}||_{L^p(\pi_{A_n})}^{q}}.
\end{equation}
\end{lemma}
\proof
If $\sup_{n\in\Nset}||1_{A_n}f_{A_n}||_{L^p(\pi_{A_n})}=\infty$, (\ref{nochungl}) is trivially satisfied.
So we may assume that $\sup_{n\in\Nset}||1_{A_n}f_{A_n}||_{L^p(\pi_{A_n})}^{q}<\infty$ .
\begin{eqnarray}
k&\le& \frac{1}{\pi(A_n)\pi(A_{n}^{c})}\int_{A_n}f_{A_n}(x)p^{2}(x,A_n^{c})\pi(dx)\nonumber\\
&\le& \frac{1}{\pi(A^{c})}\left(\int_{A_n}f_{A_n}(x)^{p}\pi_{A}(dx)\right)^{\frac{1}{p}}\left(\int_{A_n}p^{2}(x,A_n^{c})^{q}\pi_{A}(dx)\right)^{\frac{1}{q}}\nonumber\\
&\le&\frac{1}{\pi(A^{c})^{\frac{1}{p}}}\sup_{n\in\Nset}||1_{A_n}f_{A_n}||_{L^p(\pi_{A_n})}k_{2}^{\frac{1}{q}}.\nonumber
\end{eqnarray}
From here we obtain
\[
k_{2}\ge\pi(A^{c})^{\frac{q}{p}}\frac{k^{q}}{\sup_{n\in\Nset}||1_{A_n}f_{A_n}||_{L^p(\pi_{A_n})}^{q}}\ge (\frac{1}{2})^{\frac{q}{p}}\frac{k^{q}}{\sup_{n\in\Nset}||1_{A_n}f_{A_n}||_{L^p(\pi_{A_n})}^{q}}.
\]
\finishproof
Now Theorem \ref{four} follows immediately from Theorem \ref{two}, Lemma \ref{nele} and (\ref{theo4}).

As a consequence of Corollary \ref{corone} and Lemma \ref{nele} we obtain

\begin{cor}
Let $\xi_1,\xi_2,\ldots$ be a weak reversible Markov chain of order 1 with reversibility constant $C_R$ and assume that we have
that 
\begin{enumerate}
\item
$k>0$
\item
\begin{equation}\label{kapd}
\sup_{A\in\mathcal{F}:\pi(A)\ge\frac{1}{2}}\frac{p(x,A)}{p^{2}(x,A)}\le C_{\infty}<\infty\,\,\pi-a.s.
\end{equation}
\end{enumerate}
Then, for the spectrum $\sigma(P)$ of $P$ on $L_{0,1}^{2}(\pi)$, we obtain
\begin{equation}
\sigma(P)\subset B_{\sqrt{1-\frac{\kappa}{8C_{R}^{2}C_{\infty}^{2}}k^{2}}}(0).
\end{equation}  
\end{cor}
\proof
From inequality (\ref{kapd}) we obtain that from Lemma \ref{nele} with $p=\infty$, $q=1$ it follows that
$k_{2}\ge \frac{k}{C_{\infty}}$. Moreover, by Corollary \ref{corone} we see that $k_{P^{\ast}{P}}\ge \frac{k_2}{C_{R}}$.
But this yields for all $f\in L_{0,1}^{2}(\pi)$
\begin{eqnarray}
<Pf,Pf>_{\pi}&=&1-<f,(Id-P^{\ast}P)f>_{\pi}\le 1-\frac{\kappa}{8}k_{P^{\ast}{P}}^{2}\nonumber\\
&\le& 1-\frac{\kappa}{8C_{R}^{2}C_{\infty}^{2}}k^{2}\nonumber.
\end{eqnarray}
The claim follows now from the fact that $r(P)\le ||P||$, where $r(P)$ denotes the spectral radius of $P$ on $L_{0,1}^{2}(\pi)$.
\finishproof

\section{Reversibilization procedures}\label{}

We saw that in order to get information about $\sigma(P)$ of $P$, different procedures can
be used:
Choose a ``suitable'' function $h$ and consider the following expression:
\begin{equation}\label{revproc}
<f,(Id-h(P,P^{\ast}))f>_{\pi},\,\,\,f\in L_{0,1}^{2}(\pi).
\end{equation} 
Here, suitable means that $h(P,P^{\ast})$ should at least satisfy the following conditions:
\begin{itemize}
\item
$h(P,P^{\ast})$ should be self adjoint and positive
\item
\[
<f,(Id-h(P,P^{\ast}))f>_{\pi}\ge c>0,\,\,\,f\in L_{0,1}^{2}(\pi).
\]

\end{itemize}
The first condition allows a representation of (\ref{revproc}) as seen in \ref{deshalbgehtss}.
This, together with the second condition can be used to estimate the spectrum $\sigma(h(P,P^{\ast}))$ of $h(P,P^{\ast})$. Of course, $h$ must be chosen in such a way that information about $\sigma(h(P,P^{\ast}))$ 
can be used to get estimations for $\sigma(P)$. Essentially, the functions $h$ that have been 
used in the literature so far are given by
\[h_{1}(P,P^{\ast})=\frac{1}{2}(P^{\ast}+P)\,\,\,h_{2}(P,P^{\ast})=P^{\ast}P.\]
(see for example (\cite{lawler}) and (\cite{diaconis2}).
By definition, $h_1(P,P^{\ast})$ and $h_2(P,P^{\ast})$ are self adjoint. For $h_1(P,P^{\ast})$ we 
have that $<f,(Id-h_{1}(P,P^{\ast}))f>_{\pi}=\Re(<f,(Id-P)f>_{\pi})$. Since
$k=k_{h_{1}(P,P^{\ast})}$, this reversibilization yields immediate information about
the real part of $\sigma(P)$ in terms of $k$ \cite{lawler}. \\
As already seen in the weak reversible case, we use $h_{2}(P,P^{\ast})$. Compared to $h_{1}(P,P^{\ast})$, the associated isoperimetric constant has the disadvantage that it cannot be immediately related to $k$. On the other hand, we show that in some cases we are able to
compare it with $k_2$ and vice versa. Since $h_{2}(P,P^{\ast})$ is positive, this yields immediately 
in addition a global spectral gap property and not only a spectral gap property at one. \\
One may ask if there are other reasonable functions $h$ yielding sharper estimates for the spectral gap than that considered above. 
It seems to be that there are different possibilities for choosing $h$ in order to obtain  
good estimates for $\sigma(P)$. But in these cases we have to put stronger assumptions on the transition probabilities. To make things precise, let us prove Theorem \ref{fife}.

\proof
Since we assumed $P+P^{\ast}-P^{\ast}P$ to be positive, we can use the proof due to Lawler and Sokal \cite{lawler} to obtain for all $f\in L^{2}_{0,1}(\pi)$ that
\begin{equation}
<(Id-P)f,(Id-P)f>_{\pi}=<f,(Id-(P+P^{\ast}-P^{\ast}P))f>_{\pi}\ge \frac{\kappa}{8}k_{P+P^{\ast}-P^{\ast}P}^{2}.
\end{equation}
This implies that 
\[||(Id-P)f||_{2}\ge \sqrt{\frac{\kappa}{8}}k_{P+P^{\ast}-P^{\ast}P}\,\,\,\,\,\forall \,f\in L_{0,1}^{2}(\pi).\]
Moreover, we have that
\[||(Id-P)f||_2\le ||\lambda Id(f)||_{2}+||((1-\lambda)Id-P)f||_2.\]
This yields
\begin{equation}\label{ries}
||((1-\lambda)Id-P)f||_2>0\,\,\forall\,\,\lambda:\,|\lambda|<\sqrt{\frac{\kappa}{8}}k_{P+P^{\ast}-P^{\ast}P}.
\end{equation}
This implies that for all $|\lambda|<\sqrt{\frac{\kappa}{8}}k_{P+P^{\ast}-P^{\ast}P}$ we have that
$\lambda\in\sigma_{ap}(P)$.
Since we assumed that $\sigma(P)=\sigma_{ap}(P)$ we obtain 
by the spectral mapping theorem that $B_{\sqrt{\frac{\kappa}{8}}k_{P+P^{\ast}-P^{\ast}P}}$ belongs to the resolvent set of $P$. \\
It remains to show that $k_{P+P^{\ast}-P^{\ast}P}\ge k^{2}$. This can be seen as follows:
\begin{eqnarray}
k_{P+P^{\ast}-P^{\ast}P}(A)&=&\frac{1}{\pi(A)\pi(A^{c})}\left(\int_{A}p(x,A^{c})\pi(dx)+\int_{A^{c}}p(x,A)\pi(dx)\right.\nonumber\\
&&-\left.\int_{\Omega}p(x,A)p(x,A^{c})\pi(dx)\right)\nonumber\\
&=&\frac{1}{\pi(A)\pi(A^{c})}\left(\int_{A}p(x,A^{c})^2\pi(dx)+\int_{A^{c}}p(x,A)^{2}\pi(dx)\right)\nonumber\\
&\ge&\frac{1}{\pi(A^{c})}\left(\frac{1}{\pi(A)}\int_{A}p(x,A^{c})\pi(dx)\right)^{2}\nonumber\\
&&+\frac{1}{\pi(A)}\left(\frac{1}{\pi(A^{c})}\int_{A^{c}}p(x,A)\pi(dx)\right)^{2}\nonumber\\
&=&k(A)^{2}.\nonumber\\
\end{eqnarray}

\finishproof

Let us return to Proposition \ref{one} and use it to check whether the underlying Markov chain has a spectral gap or not. Let us consider 
\begin{example}
We consider the Markov chain $\xi_1,\xi_2,\ldots$ with state space $\Omega$ given by
\[\Omega=\{0\}\cup\{(a,b):a\ge 1, b\in\{1,2,\ldots,a\}\}.\]
The dynamics of the Markov chain is given by the following transition kernel:
\[p((a,b),(a,b-1))=1, \,\,\mbox{ for } b\ge2,\,\,p((a,1),0)=1,\]
$p(0,0)=\frac{1}{2}$ and 
\[p(0,(a,b))=\left\{\begin{array}{r@{\quad:\quad}l}
2^{-(a+1)}& a=b\\
0& \mbox{otherwise}\\
\end{array}\right..\]
One can show that the invariant starting distribution $\pi$ is given by
$\pi(0)=\frac{1}{2}$ and $\pi((a,b))=2^{-a+2}$ for $b\in\{1,2,\ldots,a\}$.\\
Since $p^{\ast^{n}}p^{n}((3n,2n),(3n,2n))=1\,\,\forall n\in\Nset$, it follows from Theorem \ref{one} that $P$ has not the spectral gap property.
\end{example}
\begin{remark}
The example above is due to Häggström \cite{haeg}, who used it to show that for geometrically ergodic Markov chains $\xi_1,\xi_2,\ldots$, finiteness of the second moments of a function $h$ does not ensure the validity
of a central limit theorem for the sequence  $h\circ\xi_{1},h\circ\xi_2,\ldots$. 
\end{remark}

\section{Acknowledgments}
The authors thank  Manfred Denker and Wolfgang Stadje for reading the preprint and very helpful conversations,
hints and comments.

\end{document}